\documentclass [12pt]{amsart}
\usepackage {amsmath,amssymb,amsfonts,times}
\usepackage{amscd}

\newtheorem{thm}{Theorem}[section]
\newtheorem{defin}[thm]{Definition}
\newtheorem{prop}[thm]{Proposition}
\newtheorem{lemma}[thm]{Lemma}

\newtheorem{conj}[thm]{Conjecture}

 \normalbaselines
\date{}

\newcommand \mb[1]{\mbox{$\mathbb{#1}$}}
\newcommand \mc[1]{\mbox{$\mathcal{#1}$}}

\renewcommand \P[1]{\mbox{$\mathbb{P}^{#1}$}}

\renewcommand \phi{\mbox{$\varphi$}}

\newcommand \To{\longrightarrow}

\renewcommand \rho{\mbox{$\varrho$}}
\newcommand \dd{{\bf d}}
\newcommand \ee{{\bf e}}

\newcommand \codim{\mbox{\rm codim\ }}

\def\Im{\mathop{\rm Im}\nolimits}

\title{On  Chern ratios for surfaces with ample cotangent bundle}
\author{Denis Conduch\'e and Eleonora Palmieri}

\address{Denis Conduch\'e\\
Universit\'e Louis Pasteur\\
IRMA-Mathematique\\
7, rue Ren\'e Descartes\\
67084 Strasbourg Cedex (France)
}
\email{conduche@math.u-strasbg.fr}

\address{Eleonora Palmieri\\
Dipartimento di Matematica\\
Universit\`a degli Studi Roma Tre\\
Largo San Leonardo Murialdo,1\\
00146 Roma (Italy)
}
\email{palmieri@mat.uniroma3.it}


\begin{document}
\begin{abstract}
In this paper we study the problem of density in $(1,3]$ for the Chern ratio of surfaces with ample cotangent bundle. In particular we prove density in $(1,2)$ by constructing a family of  complete intersection surfaces in a product of varieties with big cotangent bundle. We also analyse the case of complete intersections in a product of curves of genus at least 2.
\end{abstract}
\maketitle

\section{Introduction}
We work over an algebraically closed field.

The Chern numbers  of a surface of
general type satisfy the well-known inequality
$c_1^2\leq  3c_2$ due to Miyaoka (see \cite{M}). Mathematicians
tried to understand what  values the ratio $c_1^2/c_2$ can
assume.  Sommese in \cite{S} found an answer for this problem
by showing that the set of these ratios for minimal surfaces of general type is  dense in the interval $[1/5,3] $. For his proof, he uses covers of a family of
surfaces constructed by Hirzebruch in \cite{H}.

In this paper we are interested in the analogous question for surfaces with
ample cotangent bundle.
We need to recall some definitions. Given a vector bundle $E$ on a variety $X$, we define, following Grothendieck,  $\mathbb P(E)\to X$ to be the projective bundle of {\em hyperplanes} in the fibers of $E$. It is equipped with a line bundle $\mc O_{\mathbb P(E)}(1)$.

\begin{defin}
A vector bundle $E$ on a   variety $X$ is
\begin{itemize}
\item {\em ample} if the
line bundle $\mc O_{{\mathbb P}(E)}(1)$ is ample;
\item {\em big} if the
line bundle $\mc O_{{\mathbb P}(E)}(1)$ is big.
\end{itemize}
\end{defin}

Varieties with ample cotangent bundle are of general type,but they have much stronger properties: all their subvarieties are of general
type and, over the complex numbers, they are analytically
hyperbolic (any holomorphic map from $ \mathbb C$ to such a variety is constant). They 
are related to the following well-known conjecture.

\begin{conj}[{\bf Lang}]:  A smooth complex projective variety is analytically hyperbolic if and only if all its subvarieties are of general type.
\end{conj}

Surfaces with ample cotangent bundle  are minimal surfaces of general type, and they satisfy
(see \cite{FL}) the supplementary inequality $c_1^2>c_2$. So the set
of possible Chern ratios is now restricted to $(1,3]$.

As noted by  Spurr in \cite{Sp}, the surfaces constructed by Sommese do not
have, in general, this property of ampleness which  is usually quite
difficult to check.
In section \ref{se.con} we prove the following result.
\begin{thm}\label{main}
The set of ratios $c_1^2/c_2$ for surfaces with
ample cotangent bundle is dense in the interval $ (1,2)$. In other
words, any   number in $(1,2)$ is a limit of Chern ratios of surfaces with ample cotangent bundle.
\end{thm}

The proof is based on the following theorem, to be proved in section \ref{se.bog}, that generalizes a result of Bogomolov.
\begin{thm}
Let $X_1,\ldots,X_m$ be smooth projective varieties with big cotangent
bundle, all of dimension at least $d>0$. Let $Y$ be a general complete
intersection in $X_1\times\cdots\times X_m$. If $\dim
Y\leq\frac{d(m+1)+1}{2(d+1)}$,
the cotangent bundle of $Y$ is ample.
\end{thm}

In section \ref{se.linsec} we explicitly compute the Chern ratio in the special case of complete intersections  in a product of curves.

\section{Notation and conventions}

For all finite sequences of integers $\dd = (d_1,\dots,d_c)$, we set
$|\dd |=\sum_{i=1}^c d_i$. If $ V_{d_1},\dots,V_{d_c}$ are varieties, we denote by $V_\dd$ the product $V_{d_1}\times\dots \times V_{d_c}$.

We use $\sim$ for   numerical 
equivalence of divisors, and   $\cdot$ for  the intersection product of   cycles on a variety.

We denote by $c_j(\mc F)$   the $j$-th Chern class of a vector bundle $\mc F$. For any smooth variety $V$, the class $c_j(V)$ is the  $j$-th Chern class of the tangent bundle $\mc T_V$ of $V$.

\section{A generalization of Bogomolov's theorem}\label{se.bog}

We first prove a generalization of the following theorem of
Bogomolov (\cite{D}, Proposition 23).
\begin{thm}[Bogomolov]\label{thbogo}
Let $X_1,\ldots,X_m$ be smooth projective varieties with big cotangent bundle,
all of dimension at least $d>0$. Let $V$ be a general linear section of
$X_1\times\cdots\times X_m$. If $\dim V\leq \frac{d(m+1)+1}{2(d+1)}$,
the cotangent bundle of $V$ is ample.
\end{thm}

In order to generalize this result from general linear section   to
general complete intersections  we need the following.
\begin{lemma}
Let $X$ be a smooth subvariety of a projective space and let $B$ be a
subvariety of $\mathbb P(\Omega_X)$. A general complete intersection $Y$ in  $X$
of dimension at most $\frac 1 2 \codim B$ satisfies
$$
\mathbb P(\Omega_Y)\cap B=\varnothing.
$$
\end{lemma}

\begin{proof}
Let $\mathbb P^n$ be the ambient projective space and let $V_{r}=\mathbb P^{\binom{n+r}r}$ be
the Veronese variety that
parametrizes hypersurfaces of degree $r$ in $\mathbb P^n$.
Let $\dd=(d_1,\ldots,d_c)$ be a sequence of positive integers.
Now consider the variety
$$
W:=\{((t,x),Y_1,\ldots,Y_c)\in B\times V_{\dd}\,|\,x\in Y, t\in T_{X,x}\cap T_{Y_1,x}\cap\dots\cap  T_{Y_c,x}\}
$$
where $Y=X\cap Y_1\cap \dots \cap Y_c$. 
We have   two projections $W\to B$ and  $W\to V_{\dd}$. The fibers of
 the first projection have
codimension $2c$, so the second projection is not dominant
 when $2c>\dim B$. 
But $c=\codim_X Y=\dim X-\dim Y$, so the last condition is equivalent to
$$
2(\dim X-\dim Y)-1\geq \dim B=\dim \mathbb P(\Omega_X)-\codim B=2\dim
X-1-\codim B
$$
which proves the lemma.
\end{proof}

Analogously we can prove the following.
\begin{lemma}
Let $Y$ be a general complete intersection in a product $X_1\times X_2$ in
a projective space. If $2\dim Y\leq \dim X_1+1$, the projection
$Y\to X_1$ is finite.
\end{lemma}

\begin{proof}
We consider in the ambient space $\mathbb P^n$ the closure of the
locus
$$
\{(x_1,x_2,x_2',Y_1,\ldots, Y_c)\in X_1\times X_2\times X_2\times
V_{\dd}\,|\,x_2\neq x_2',\,(x_1,x_2),(x_1,x_2')\in
Y_1\cap\dots\cap Y_c\}
$$
This variety has a natural projection to $X_1\times X_2\times X_2$ whose
fibers have codimension $2c$. As in the previous lemma the projection to
$V_{\dd}$ has fibers of dimension at most
1 if the condition $2c\geq \dim (X_1\times X_2\times X_2)-1$ is satisfied. When
$c=\codim Y$, this amounts to saying that the projection $Y\to X_1$ is
finite and that
$$
2\dim (X_1\times X_2)-2\dim Y\geq \dim (X_1\times X_2\times X_2)-1
$$
which is equivalent to
$$
\dim Y\leq \dim X_1+1
$$
as wanted.
\end{proof}

We now obtain the desired statement.
\begin{thm}\label{bog2}
Let $X_1,\ldots,X_m$ be smooth projective varieties with big cotangent
bundle, all of dimension at least $d>0$. Let $Y$ be a general complete
intersection in $X_1\times\cdots\times X_m$. If $\dim
Y\leq\frac{d(m+1)+1}{2(d+1)}$,
the cotangent bundle of $Y$ is ample.
\end{thm}

\begin{proof}
The proof goes as in Bogomolov's theorem (see \cite{D}).
\end{proof}

\begin{proof}[Remark]
Of course, in case the embedding of $X_1\times\cdots\times X_m$ in a projective space comes from embeddings of each $X_i$ in a projective space followed by a Segre embedding, this theorem is just a consequence of Bogomolov's theorem \ref{thbogo}.
\end{proof}

\section{Computation of Chern numbers}\label{se.con}

We compute the Chern numbers of a complete intersection surface in a product.

\begin{prop}\label{propc}
Let $X_1,\ldots,X_m$ be smooth projective varieties with big cotangent
bundle, all of dimension at least $d>0$, such that $2
 \leq\frac{d(m+1)+1}{2(d+1)}$. Let $N$ be the dimension of $X=X_1\times\cdots\times X_m$. Embed $X$ in a projective space $\mathbb P^n$.
Let $S$ be a surface, general complete intersection of degree
 $\dd=(d_1,\dots,d_{N-2})$ in $X$.
We   have 
$$
\begin{array}{rcl}
c_1^2(S)&=& c_1^2(X) \vert_S -2 |\dd| c_1(X)\vert_S \cdot H \vert_S + |\dd|^2 H^2\vert_S\\
c_2(S)  &=& c_2(X) \vert_S - |\dd| c_1(X)\vert_S \cdot H \vert_S +
 \left(\sum_{i\leq j}d_i d_j\right) H^2 \vert_S
\end{array}
$$
where $H$ is a   hyperplane in $\mathbb P^n$.
\end{prop}
\begin{proof}

According to Theorem \ref{bog2}, the cotangent bundle of $S$ is ample.
The surface $S$ is the intersection of $X$ with  hypersurfaces
$L_1,\ldots,L_{N-2}$ of degrees $d_1,\ldots,d_{N-2}$ in  $\mathbb P^n$. We want to compute the ratio $c_1^2/c_2$ for
$S$. Since $S$ is smooth, we have a short exact sequence of sheaves (see for example \cite[p.182]{Ha})
$$
0\To \mc T_S \To (\mc T_X)\vert_S \To \mc N_{S/X} \To 0
$$ 
and we also know that, for a complete intersection,
$$
\mc N_{S/X} \simeq \left(\mc I_S/\mc I_S^2\right)^*\simeq \bigoplus_{i=1}^{N-2} \mc
O_S(L_i)
$$ 
where   $\mc I_S$ is  the ideal sheaf of   $S$ in $X$.
So, for the Chern classes,
$$
 c(\mc T_X)\vert_S =c(\mc T_{S})\cdot c\left(\mc N_{S/X}\right)
$$
$$
c\left(\mc N_{S/X}\right)=c\biggl(\bigoplus_{i=1}^{N-2} \mc
O_S(L_i)\biggr) = \prod_{i=1}^{N-2}(1+c_1(\mc O_S(L_i)))
$$
Let us compute $c_1^2(S)$:
$$
c_1(\mc T_X)\vert_S=c_1(\mc T_S)+\sum_{i=1}^{N-2}c_1(\mc O_S(L_i))
$$
which implies, with $c_1(S)=c_1(\mc T_S)$ and $c_1(X)=c_1(\mc
T_{X})$,
$$ 
c_1(S) =  c_1(X)\vert_S -\sum_{i=1}^{N-2}c_1(\mc O_S(L_i))
$$
Now we remember that, if $H$ is a hyperplane in $\mathbb P^n$,
we have for all $i\in\{1,\ldots,N-2\}$, $L_i\sim d_iH$, so that 
\begin{eqnarray*} 
c_1(S)&=&c_1(X)\vert_S-|\dd| H\vert_S 
\end{eqnarray*}
and
\begin{eqnarray*} 
c_1^2(S)&=&(c_1(X)-|\dd| H\vert_X)\vert_S^2\\
&=&c_1^2(X)\vert_S-2 |\dd| c_1(X)\vert_S\cdot H\vert_S+|\dd|^2H^2\vert_S 
\end{eqnarray*}
Let us compute $c_2(S)$:
$$
c_2(\mc T_X)\vert_S=c_2(\mc T_S)+c_1(\mc T_S)\cdot\sum_{i=1}^{N-2}c_1(\mc
O_S(L_i))+\sum_{i<j}c_1(\mc O_S(L_i))\cdot c_1(\mc
O_S(L_j))
$$
so that
\begin{eqnarray*} 
c_2(S)&=&c_2(\mc T_X)\vert_S-c_1(\mc T_S)\cdot \sum_{i=1}^{N-2}c_1(\mc
 O_S(L_i))-\sum_{i<j}c_1(\mc O_S(L_i))\cdot c_1(\mc O_S(L_j))\\
&=&c_2(X)\vert_S-c_1(S)\cdot\sum_{i=1}^{N-2}L_i\vert_S-\sum_{i<j}(L_i\cdot L_j)\vert_S\\
&=&\Bigl(c_2(X)-(c_1(X)-|\dd| H\vert_X)\cdot\sum_{i=1}^{N-2}L_i\vert_X-\sum_{i<j}(L_i\cdot L_j)\vert_X\Bigr)\vert_S\\
&=&\Bigl(c_2(X)-|\dd|(c_1(X)-|\dd| H\vert_X)\cdot  H\vert_X- \sum_{i<j}d_i
 d_j H^2\vert_X\Bigr)\vert_S
\end{eqnarray*}
So, if we set
\begin{equation}\label{ab}
a:=c_1(X)\cdot H^{N-1}\vert_X
\qquad\hbox{and}\qquad
b:=H^N\vert_X
\end{equation}
we can write
$$
c_1^2(S)= \left(\prod_{i=1}^{N-2}d_i\right)\left( c_1^2(X) \cdot
 H^{N-2}\vert_X -2a|\dd|+b|\dd|^2 \right)
$$
and
\begin{eqnarray*}
c_2(S)&=&\left(\prod_{i=1}^{N-2}d_i\right)\left(c_2(X)\cdot
H^{N-2}\vert_X-a|\dd|+b\Big(|\dd|^2-\sum_{i<j}d_i d_j\Big) \right)\\
&=&\left(\prod_{i=1}^{N-2}d_i\right)\left(c_2(X)\cdot H^{N-2}\vert_X-a|\dd|+b\sum_{i\leq j}d_id_j \right)
\end{eqnarray*}
\end{proof}

Assume moreover that the $d_i$  are multiple of the same integer
$d$. So $d_i=e_id$ for all $i\in\{1,\dots,N-2\}$. The ratio is
\begin{eqnarray}\label{eq.ratio}
\frac{c_1^2(S)}{c_2(S)}&=&\frac{c_1^2(X)\cdot
H^{N-2}\vert_X-2a|\dd|+b|\dd|^2}{c_2(X)\cdot H^{N-2}\vert_X-a|\dd|+b\sum_{i\leq j}d_i
d_j}\\
&=&\frac{c_1^2(X)\cdot
H^{N-2}\vert_X-2ad|\ee|+bd^2|\ee|^2}{c_2(X)\cdot
H^{N-2}\vert_X-ad|\ee|+bd^2\sum_{i\leq j}e_i
e_j}\notag
\end{eqnarray}
Letting $d$ go to infinity, we obtain
$$
\lim_{d\to+\infty}\frac{c_1^2(S)}{c_2(S)}=\frac{|\ee|^2}{\displaystyle\sum_{1\le i\leq j\le N-2}e_i
e_j}
$$
which is   a rational number between $1$ and $2$.

\begin{thm}\label{valthm}
The set of values of the fraction above is dense in the interval $(1,2)$.
\end{thm}

We want to rephrase Theorem \ref{valthm} as follows. Set $T:=|\ee|$ and $U:= e_1^2+\dots+e_{N-2}^2$. Then 
$$
\frac{|\ee|^2}{\sum_{1\le i\leq j\le N-2} e_ie_j}=\frac{T^2}{(T^2+U)/2}=\frac{2T^2}{T^2+U}
$$
The statement of the theorem is equivalent to ask for the values of
$$
\frac{T^2+U}{2T^2}=\frac12\Bigl(1+\frac U{T^2}\Bigr)
$$
to be dense in $(1/2,1)$ or for the values of
$$
\frac{T^2}U
$$
to be dense in $(1,+\infty)$.
 The new formulation is the following.
\begin{thm} \label{valthm1}
The set of rational numbers
$$
\left\{ \frac{\Bigl(\sum_{i=1}^M e_i\Bigr)^2}{\sum_{i=1}^M e_i^2}\ \Bigg|\
 M\ge4,\ e_1,\dots,e_M\in\mb N^*\right\}
$$
is dense in the interval $(1,+\infty)$.
\end{thm}

\begin{proof}
Fix $M\ge4$ and define a function 
$$
\begin{matrix}
 f_M:&(\mb  R^{+*})^M&\To &\mb R^+\\
&(e_1,\dots,e_M)&\longmapsto &\frac{\Bigl(\displaystyle\sum_{i=1}^M e_i\Bigr)^2}{\displaystyle\sum_{ i=1}^{M}e_i^2}
\end{matrix}
$$
Observe that for all real numbers $\lambda>0$, 
$$
 f_M(\lambda e_i,\dots,\lambda e_M)= f_M(e_1,\dots,e_M)
$$
This implies   $ f_M( (\mb  N^*)^M) =f_M ((\mb  Q^{+*})^M)$. Since $\mb Q^+$ is dense in $\mb R^+$, $   f_M((\mb  Q^{+*})^M )$ is dense in $  f_M( (\mb  R^{+*})^M) $.
To prove  the theorem, we have only to show that 
$$
\bigcup_{M\ge4} \Im( f_M)=(1,+\infty) 
$$
It is easy to see that $  \Im( f_M)\subseteq (1,+\infty)$. On the other hand, since we have 
$$
 f_M(1,\dots,1)=\frac{M^2}{M}=M
$$ 
and 
\begin{eqnarray*}
  f_M(1,\dots,1,M^2)&=&\frac{(M-1+M^2)^2}{M-1+M^4}=\frac{M^4+M^2+1-2M+2M^3-2M^2}{M^4+M-1}\\
&=&1+\frac{2M^3-M^2-3M+2}{M^4+M-1}=:1+\varepsilon(M)
\end{eqnarray*}
the whole interval $[1+\varepsilon(M),M]$ is contained in the image of $  f_M$. So
$$
\bigcup_{M\ge4} \Im( f_M)\supseteq\bigcup_{M\ge4}[1+\varepsilon(M),M]=(1,\infty)
$$
and the theorem is proved.
\end{proof}

Note that the hypothesis $2\leq\frac{d(m+1)+1}{2(d+1)}$ in Proposition \ref{propc} is equivalent to
 $$
m\ge 3+\frac3d
 $$
 Taking for the $X_i$ curves of genus at least $2$ and $m=N\ge 6$, we see that Theorem \ref{main} is proved.
\vskip .5 cm

We also have, as a corollary, the following result.
\begin{prop}\label{finite}
For any fixed product of projective varieties $X=X_1\times\dots\times X_m$ as above, there is only a finite number of ratios $\frac{c_1^2(S)}{c_2(S)}$ greater than or equal to $2$ for a general complete intersection surface  $S$ in $X$.
\end{prop}
\vskip .3 cm

We recall the following general fact (see, for example, \cite{BPV}):
\begin{prop}
If $F$ is any surface of general type,   $c_2(F)>0$. 
\end{prop}
 \vskip .5 cm

\begin{proof}[Proof of Proposition \ref{finite}.]
With the notation of proposition \ref{propc}, we have
$$
\frac{c_1^2(S)}{c_2(S)}=\frac{c_1^2(X)\cdot
H^{N-2}\vert_X-2a|\dd|+b|\dd|^2}{c_2(X)\cdot H^{N-2}\vert_X-a|\dd|+b\sum_{i\leq j}d_i
d_j}
$$
so $c_1^2(S)\ge 2c_2(S)$ if and only if 

$$
c_1^2(X)\cdot H^{N-2}\vert_X-2a|\dd|+b|\dd|^2\geq 2\Bigl(c_2(X)\cdot H^{N-2}\vert_X-a|\dd|+b\sum_{i\leq j}d_i
d_j\Bigr)
$$
or, equivalently, if and only if
$$
\Bigl(c_1^2(X)-2c_2(X)\Bigr)\cdot
 H^{N-2}\vert_X\geq b\Bigl(2\sum_{i\leq j}d_i
d_j-|\dd|^2\Bigr)=b\sum_{i=1}^{N-2}d_i^2=\Bigl(\sum_{i=1}^{N-2}d_i^2\Bigr) H^N\vert_X
$$
Since $H^N\vert_X>0$, we have
$$
\sum_{i=1}^{N-2}d_i^2\le\frac{\Bigl(c_1^2(X)-2c_2(X)\Bigr)\cdot
H^{N-2}\vert_X}{H^N\vert_X}
$$
We have now to consider two cases:
\begin{enumerate}
\item[$\bullet$] if $(c_1^2(X)-2c_2(X))\cdot H^{N-2}\vert_X\le 0$,
  the set of possible vectors $\dd$ is empty;
\item[$\bullet$] if $(c_1^2(X)-2c_2(X))\cdot H^{N-2}\vert_X=\alpha>0$, every $d_i$ is bounded by 
$$
 \sqrt{\frac{\alpha}{H^N\vert_X}} 
$$
so the number of vectors $\dd$ is finite.
\end{enumerate}
\end{proof}

\section{Explicit computations for linear sections of a product of curves}\label{se.linsec}

Let us fix   curves $X_1,\dots,X_N$ of respective genera $g_i\geq 2$ and embed each  $X_i$    in a projective space via a multiple  $l_iK_{X_i}$ of the canonical bundle, with $l_i\geq 1$ ($\geq 2$ if $X_i$ is hyperelliptic) as a curve of degree $\alpha_i=l_i(2g_i-2)$. 

Set $X=X_1\times\dots\times X_N$ and embed $X$ in a projective space $\mathbb P^n$ via a Segre embedding. If $\pi_i:X\to X_i$ is the projection on the $i$-th factor and $K_{\widehat i} =\pi_i^*K_{X_i}$, and if $H$ is a hyperplane in $\mathbb P^n$, we have
$$
H\vert_X=\sum_{i=1}^N l_iK_{\widehat i}
$$
For each $i\in\{1,\ldots,N\}$,  let $p_i$ be a point in $X_i$. For all $j\in\{1,\dots,N\}$  
and for all multi-indices
$$
I_j=(i_1,\ldots,i_j)
$$
 with $1\le i_1<\ldots<i_j\le N$, set
$$
X_{I_j}:=X_{i_1}\times\cdots\times X_{i_j}
$$
With this notation, we have for each $I_j$  a projection
$$
\pi_{I_j}:X \To X_{I_j}
$$
Let $X_{\widehat{I_j}}$ be the fibre of $\pi_{I_j}$ over the point $(p_{i_1},\ldots,p_{i_j})$.
We have
$$
H\vert_X= \sum_{i=1}^N \alpha_{i} X_{\widehat{(i)}}
$$
We need to compute $H^j\vert_X$:
$$
H^j\vert_X=j!\sum_{I_j}\Big(\prod_{k\in I_j}\alpha_k\Big) X_{\widehat{I_j}}
$$
Writing $X_i\times X_j $ instead of $X_{\widehat{\{1,\dots,N\}\setminus\{i,j\}}}$, we obtain
\begin{eqnarray*}
H^{N}\vert_X&=&N!\prod_{i}\alpha_i
\\
H^{N-1}\vert_X&=&(N-1)!\sum_{i}\Bigl(\prod_{k\neq i}\alpha_k\Bigr)X_i
\\
H^{N-2}\vert_X&=&(N-2)!\sum_{i<j}\Bigl(\prod_{k\neq i,j}\alpha_k\Bigr)(X_i\times X_j)
\end{eqnarray*}
Moreover,  
$$
c_1(X)=\sum_i \pi_i^*c_1(X_i)=-\sum_i K_{\widehat i}\sim -\sum_i (2g_i-2)X_{\widehat{(i)}}
$$
and
$$
c_2(X)=\sum_{i<j}\pi_i^*c_1(X_i)\cdot\pi_j^*c_1(X_j)=\sum_{i<j} 
K_{\widehat i}\cdot K_{\widehat j}
\sim \sum_{i<j} (2g_i-2)(2g_j-2)X_{\widehat{(i,j)}}
$$
Since $c_1^2(X_i)=0$ for all $ i\in\{1,\ldots,N\}$, we always have $c_1^2(X)=2c_2(X)$ in $H^4(X,\mathbb Z)$.

\begin{prop}
A  general complete intersection surface $S$ in a product of curves  satisfies the inequality
$$
c_1^2(S)<2c_2(S)
$$
\end{prop}

\begin{proof}
It is just a consequence of Proposition \ref{finite} and of $c_1^2(X)=2c_2(X)$ for $X$ a product of curves.
\end{proof}
\vskip .3 cm

In this case we can compute Chern classes more explicitly and obtain a precise numerical result.

\begin{prop}\label{pr.ratio}
Let $S$ be a
complete intersection surface of degree $\dd=(d_1,\ldots,d_{N-2})$ in a product of curves $X=X_1\times\dots\times X_N$ such that $X_i$ is embedded in the projective space via the $l_i$-canonical bundle. Then 
$$
\frac{c_1^2(S)}{c_2(S)}=2-\frac{N(N-1)\Bigl(\sum_{i=1}^{N-2}d_i^2\Bigr)}
{\sum_{1\le i<j\le N}\frac 1{l_il_j}+(N-1)\Bigl(|\dd|\sum_{i=1}^{N}\frac 1{l_i}+N \sum_{1\le i\leq j\le N-2}d_id_j\Bigr)}
$$
In particular, when $S$ is a linear section, we have
$$
\frac{c_1^2(S)}{c_2(S)}=2-\frac{N(N-1)(N-2)}
{\sum_{1\le i<j\le N}\frac1{l_il_j}+(N-1)(N-2)\Bigl(\frac{N(N-1)}{2}+\sum_{i=1}^N\frac 1{l_i}\Bigr)}
$$
\end{prop}

Recall that we need to take $N\ge 6$ in order to be able to apply  Bogomolov's theorem to get a surface $S$ with ample cotangent bundle.

\begin{proof}
We prove the first assertion. The second one follows when we set $\dd=(1,\ldots,1)$.
Assume $S$ is a complete intersection in $X$, so that $S\sim \left(\prod_{i=1}^{N-2} d_i\right)H^{N-2}\vert_X$.
We have, using notation (\ref{ab}),
\begin{eqnarray*}
a&=&c_1(X)\cdot H^{N-1}\vert_X \\
&=&-(N-1)!\Bigl( \sum_{i=1}^N (2g_i-2) X_{\widehat{(i)}}\Bigr)\cdot \Bigl(\sum_j\Bigl(\prod_{k\neq j}\alpha_k\Bigr) X_j\Bigr)\\
&=&-(N-1)!\sum_j\Bigl((2g_j-2)\prod_{k\neq j}\alpha_k\Bigr)\\
&=&-(N-1)!\Bigl(\prod_{i}(2g_i-2)\Bigr)\sum_j\prod_{k\neq j}l_k
\end{eqnarray*}
and
\begin{eqnarray*}
b=H^N\vert_X=N!
\prod_{i}\alpha_i=N!
 \Bigl(\prod_{i=1}^{N}(2g_i-2)\Bigr)  \Bigl (\prod_{i=1}^{N}l_i\Bigr)
\end{eqnarray*}

Moreover,
\allowdisplaybreaks{
\begin{eqnarray*}
c_2(X)\cdot H^{N-2}\vert_X&=&( \sum_{i<j} (2g_i-2)(2g_j-2)X_{\widehat{(i,j)}})\cdot H^{N-2}\vert_X\\
&=&(N-2)!
\sum_{i<j} (2g_i-2)(2g_j-2)\Bigl(\prod_{k\ne i,j}\alpha_k\Bigr)\\
&=&(N-2)!
\Bigl(\prod_{i=1}^N (2g_i-2)\Bigr)\sum_{i<j}\Bigl(\prod_{k\ne i,j}l_k\Bigr)
\end{eqnarray*}
}

Now it follows from Proposition \ref{propc} and the fact that $c_1^2(X) =2c_2(X) $ that we have
\begin{eqnarray*}
\frac{c_1^2(S)}{c_2(S)}&=&2+\frac{b(|\mathbf d|^2-2\sum_{i\leq j}d_id_j)}{c_2(X)\cdot H^{N-2}\vert_X-a|\mathbf d|+b\sum_{i\leq j}d_id_j}\\
&=&2+\frac{-N(N-1)\left(\sum_i d_i^2\right)\prod_i l_i}{\sum_{i<j}\prod_{k\neq i,j}l_k+(N-1)\Bigl(|\mathbf d|\sum_i\prod_{j\neq i} l_j+N
\left(\sum_{i\leq j}d_id_j\right)\prod_i l_i\Bigr)}\\
&=&2-\frac{N(N-1)(\sum_{i}d_i^2)}
{\sum_{i<j}\frac 1{l_il_j}+(N-1)\Bigl(|\dd|\sum_i\frac 1{l_i}+N \sum_{i\leq j}d_id_j \Bigr)}
\end{eqnarray*}
where in the second equality we simply observe that all terms are multiples of 
$$
(N-2)!\Bigl(\prod_{i=1}^N(2g_i-2)\Bigr)
$$
and in the third equality we divide both the numerator and the denominator by $\Bigl(\prod_{i=1}^{N-2}l_i\Bigr)$.
\end{proof}

\section*{Remarks and open problems}

It is interesting to note that, to the best of our knowledge, not many examples of surfaces with ample cotangent bundle are known. Some of them come from Hirzebruch's construction, which uses arrangements of lines in \P 2, and satisfy a criterion for ampleness which can be found in \cite{S}. Another such criterion is proved by Spurr in \cite{Sp} to construct double covers of some special Hirzebruch surfaces.
For most of the above examples the Chern ratio is not greater than 2, and we know only a few sporadic examples with $c_1^2\ge 2c_2$, while a result of Miyaoka (see \cite{M2}) shows the ampleness of the cotangent bundle for  surfaces with $c_1^2=3c_2$, i.e., quotients of the unit ball.

So the question about density in $[2,3]$ is still open.

\vskip .5 cm

\section*{Acknowledgments}
We wish to thank Professors Olivier Debarre and Lucia Caporaso for having introduced us to this subject, and in particular Professor Debarre for his guidance during the preparation of this paper. We also wish to thank Professor Alfio Ragusa and the other organizers of the EAGER summer school PRAGMATIC 2004 for the  nice stay in Catania where this work started. 

Moreover both authors would like to heartily thank Ren\'{e} Schoof for
many useful discussions and remarks, and for having suggested the proof of Theorem \ref{valthm1}.

\end{document}